\documentclass[12pt]{amsart}
\usepackage{mathrsfs}
\usepackage[colorlinks]{hyperref}
\usepackage{xcolor}

\usepackage{setspace}

\usepackage{mathscinet}
\usepackage{latexsym}
\usepackage{amsthm}
\usepackage{amssymb}
\usepackage{amsfonts}
\usepackage{amsmath}
\usepackage{verbatim}

\newtheorem{theorem}{Theorem}[section]
\newtheorem{lemma}[theorem]{Lemma}
\newtheorem{proposition}[theorem]{Proposition}
\newtheorem{corollary}[theorem]{Corollary}
\newtheorem{hypothesis}[theorem]{Hypothesis}
\theoremstyle{definition}
\newtheorem{definition}[theorem]{Definition}
\newtheorem{example}[theorem]{Example}
\theoremstyle{remark}
\newtheorem{remark}[theorem]{Remark}
\numberwithin{equation}{section}
\def\d{{\rm d}}
\newcommand{\e}{{\rm e}}

\begin{document}

\title[Differentiability of transition semigroup]{Differentiability of transition semigroup of generalized Ornstein-Uhlenbeck process: a probabilistic approach}

\author{Ben Goldys}
\address{School of Mathematics and Statistics, The University of Sydney, Sydney 2006, Australia}
\email{Beniamin.Goldys@sydney.edu.au}

\author{Szymon Peszat}
\address{Szymon Peszat,  Institute of Mathematics, Jagiellonian University, {\L}ojasiewicza 6, 30--348 Krak\'ow, Poland}
\email{napeszat@cyf-kr.edu.pl}


\begin{abstract}
Let $P_s\phi(x)=\mathbb{E}\, \phi(X^x(s))$,  be the  transition semigroup on the space $B_b(E)$ of bounded measurable functions on a Banach space $E$, of the Markov family  defined by the linear equation with additive noise 
$$
\d X(s)= \left(AX(s) + a\right)\d s + B\d W(s), \qquad X(0)=x\in E. 
$$
We give a simple probabilistic proof of the fact that null-controlla\-bility of the corresponding deterministic system 
$$
\d Y(s)= \left(AY(s)+ B\mathcal{U}(t)x)(s)\right)\d s, \qquad Y(0)=x,
$$
implies that for any $\phi\in B_b(E)$, $P_t\phi$ is infinitely many times Fr\'echet differentiable and that 
$$
D^nP_t\phi(x)[y_1,\ldots ,y_n]= \mathbb{E}\, \phi(X^x(t))(-1)^nI^n_t(y_1,\ldots, y_n),
$$
where $I^n_t(y_1,\ldots,y_n)$ is the symmetric n-fold It\^o integral of the controls $\mathcal{U}(t)y_1,\ldots \mathcal{U}(t)y_n$. 
\end{abstract}
\thanks{This  work was 
supported by  the Australian Research Council Project DP200101866. 
} 

\date{}
\maketitle
\tableofcontents
\section{Introduction}
We will consider the following linear stochastic differential equation 
\begin{equation}\label{eq_intro1}
\d X(s)= \left(AX(s) + a\right)\d s + B\d W(s), \qquad X(0)=x\in E\,,
\end{equation}
where $E$ is a Banach space and $A$ generates a strongly continuous semigroup on $E$ (see Section \ref{sec1} for precise formulation). Under conditions specified in Section \ref{sec1} this equation has, for every $x\in E$, a unique mild solution $X^x$  known as a (generalised) Ornstein-Uhlenbeck process. Since the family $\left\{X^x;\,\,x\in E\right\}$ is Markov in $E$, we can define the corresponding transition semigroup $P_t\phi(x)=\mathbb E\,\phi(X^x(t))$. Investigation of differentiability of $P_t\phi$ was initiated by Kolmogorov in \cite{kolm1}, where a simple version of equation \eqref{eq_intro1} is considered in $\mathbb R^2$. If $E=\mathbb R^d$ and the pair $(A,B)$ satisifies the H\"ormander hipoellipticity condition (see \cite{hor}) then the function $P_t\phi$ is Fr\'echet differentiable for every bounded Borel function $\phi$ and 
\begin{equation}\label{eq_intro2}
\|DP_t\phi\|_\infty\le\frac{c}{t^{k/2}}\|\phi\|_\infty,\quad t\le 1\,,
\end{equation}
where the value of $k$ follows from the H\"ormander condition as well. 
\par
It has been known for a long time that the H\"ormander condition holds if and only if the deterministic controlled system 
\[\frac{dY}{dt}(t)=AY(t)+Bu(t),\quad Y(0)=x\]
is null-controllable for every $x\in\mathbb R^d$ and if and only if the transition semigroup $(P_t)$ is strongly Feller. It was proved in \cite{dzsf} that this formulation can be carried to Hilbert spaces,  where the estimate \eqref{eq_intro2} takes the form 
\begin{equation}\label{eq_intro3}
\||DP_t\phi\|_\infty\le\|\mathcal U(t)\|\|\phi\|_\infty,
\end{equation}
with $\|\mathcal U(t)\|$ being the norm of the controllability operator.

 In many problems it is important to estimate the gradient of the transition semigroup in  function spaces such as weighted $L^p$ spaces or weighted Sobolev and Besov spaces, in which case \eqref{eq_intro3} is not satisfactory. The aim of this paper is to derive a pointwise formula for $DP_t\phi$ that may be used to derive norm estimates in various function norms. Our main result, see Theorem \ref{T17}, provides a probabilistic formula 
 \begin{equation}\label{eq_point}
D^nP_t\phi(x)[y_1,\ldots y_n]= \mathbb{E}\, \phi(X^x(t))(-1)^nI^n_t(y_1,\ldots, y_n),
\end{equation}
where $I^n_t(y_1,\ldots,y_n)$ is the symmetric n-fold It\^o integral of the controls $\mathcal{U}(t)y_1,\ldots \mathcal{U}(t)y_n$ driving points $y_1,\ldots, y_n$ in $E$ to zero. This result is obtained under rather mild assumptions on $(A,B)$. In particular, the process $BW$ (hence the null control $Bu$) can take values in a Banach space bigger than $E$ (see Section \ref{sec1} for details) if the $C_0$-semigroup generated by $A$ has appropriate smoothing properties. As a result, we obtain smoothing properties of $P_t$ and the estimate \eqref{eq_intro3} as well as the estimates for $\|D^nP_t\phi \|_\infty $. One of the main results of the paper are the  formulae and estimates for $D^nP_t\phi$ in the case of boundary noise studied recently in \cite{Goldys-Peszat}. Let us note that the case of boundary noise requires use of weighted $L^p$-spaces, where the pointwise formula \eqref{eq_point} for the gradient of the transition semigroup has to be used. 

We easily recover the known results \eqref{eq_intro2}, \eqref{eq_intro3} in the hypoelliptic case in $\mathbb R^d$. Theorem \ref{T17} will also allow us to give rapid proofs of some known facts such as the Liouville property and absolute continuity of transition measures. 

While some of our results are known, the method of proof based on the Girsanov theorem seems to be new. Let us recall that a probabilistic formula for the gradient of the transition semigroup is also provided by the famous Bismut--Elworthy--Li formula \cite{bel} extended later to Hilbert spaces in \cite{Peszat-Zabczyk}. However, this formula requires the operator $B$ in \eqref{eq_intro1} to be invertible for $E=\mathbb R^d$ or the process $BW$ close to cylindrical Wiener process if $E$ is a Hilbert space. We do not require these assumptions.

Additionally, we will provide detailed arguments for the null controllability of equation \eqref{eq_intro2} in the case of boundary control. 
\section{Formulation and main results}\label{sec1}
Let $(W_k)$ be a finite or infinite sequence of independent standard real-valued Wiener processes defined on a probability space $(\Omega,\mathfrak{F},\mathbb{P})$. Let $(E, \|\cdot \|_{E})$ and $(V,\|\cdot\|_V)$ be Banach spaces. We assume that $E$ is continuously and densely embedded into $V$. Let $(H, \langle \cdot,\cdot \rangle_H)$, be a Hilbert space endowed with an orthonormal basis $(e_k)$. We will consider the following linear stochastic evolution equation with an additive noise 
\begin{equation}\label{E11}
\d X(s)= \left(AX(s)+ a\right)\d s +   B\d W(s), \qquad X(0)=x\,.
\end{equation}
We assume that 
$$
W(s)= \sum _k W_k(s)e_k, 
$$
is a cylindrical Wiener process on $H$, 
  $B$ is a linear operator from $H$ into $V$, $(A,\text{Dom}(A))$ is the  generator of a $C_0$-semigroup $S=(S(s); \, s\ge 0)$ on $E$, and $a\in V$. 

Our first assumption ensures  the existence of a solution to \eqref{E11} in $E$.  
\begin{hypothesis}\label{H11}
For each $s>0$, $S(s)$ has  an extension to a bounded linear operator from $V$ to $E$.  Moreover, for each $k$ and $s>0$, the mapping 
$$
[0,s)\ni r\mapsto S(s-r)Be_k 
$$
is integrable in $E$ with respect to $W_k$, and the series 
$$
\sum_k \int_0^s S(s-r)Be_k \d W_k(r)
$$
converges in probability in $E$. We also assume that for any $s>0$ the integral  $\int_0^s S(r) a\d r$ converges in $E$.
\end{hypothesis}

\begin{remark}\label{R12}
{\rm If $E$ is a Hilbert space then the first part of Hypothesis \ref{H11} can be formulated equivalently: for each $r>0$,  $S(r)B$ is a Hilbert--Schmidt operator from $H$ into $E$ and 
$$
\int_0^s \|S(r)B\|^2_{L_{(HS)}(H,E)}\d r <+\infty, \qquad s>0,
$$
where $\|\cdot \|^2_{L_{(HS)}(H,E)}$ stands for the Hilbert--Schmidt operator norm. Explicit conditions for stochastic integrability in $L^p$-spaces are also known, see for example  Proposition A.1 from \cite{Brzezniak-Veraar}.  
}
\end{remark}
Clearly, under Hypothesis \ref{H11}, equation \eqref{E11} has a unique mild solution 
$$
X^x(s)= S(s)x +\int_0^s S(s-r)a\d r+  \sum_k \int_0^t S(s-r)Be_k \d W_k(r). 
$$
In particular 
$$
X^0(s)=  \int_0^t S(s-r)a\d r+ \sum_k \int_0^t S(s-r)Be_k \d W_k(r)
$$
is a solution starting from $x=0\in E$. 

Moreover, \eqref{E11} defines on the space $B_b(E)$ of bounded measurable functions on $E$, the transition semigroup 
$P_s\phi(x)= \mathbb{E}\, \phi\left(X^x(s)\right)$, $\phi\in B_b(E)$, $s\ge 0$, $x\in E$.  
\subsection{Some definitions}

Let $L^2(0,T;H)$ denote the Hilbert space of measurable  square integrable mappings $u\colon [0,T]\mapsto H$. 

Consider the following   deterministic system 
\begin{equation}\label{E12}
\d Y(s)= \left(AY(s) +  Bu(s)\right)\d s,  \qquad Y(0)=x. 
\end{equation}
\begin{definition}
The deterministic system \eqref{E12}  is said to be \emph{null-controllable at time $T$} if  for any $t\ge T$ there exists a bounded linear operator $\mathcal{U}(t)\colon E\mapsto L^2(0,t;H)$,  such that:
$$
S(t)x + \sum_k \int_0^ t S(t-s)Bu(t,x; s)\d s=0,
$$
where $u(t,x;s)= \left(\mathcal{U}(t)x \right)(s)$. 
\end{definition}

It is easy to see that the system \eqref{E12} is null controllable if and only if the system 
\[\d Y(s)= \left(AY(s) +a+  Bu(s)\right)\d s,  \qquad Y(0)=x\]
is null-controllable.

Connections between null-controllability and strong Feller property of linear stochastic evolution equations on Hilbert spaces with additive noise were establish in the following  series of papers and books: Zabczyk \cite{Zabczyk}, Da Prato and Zabczyk \cite{DaPrato-Zabczyk2, DaPrato-Zabczyk3}, and Da Prato, Pritchard, and Zabczyk \cite{DaPrato-Pritchard-Zabczyk}. 


\begin{definition}
{\rm System \eqref{E12} is \emph{null-controllable with uniform vanishing energy}, if it is null controllable at some time $T$, and the  controllability operator $\mathcal{U}(t)\to 0$ in the operator norm as $t\to +\infty$. It is \emph{null-controllable with vanishing energy} if $\mathcal{U}(t)x\to 0$ in the $L^2(0,t; H)$ norm as $t\to +\infty$ for any $x\in E$.  }
\end{definition}

The connections between null controllability with vanishing energy and the Liouville  type property  was established in \cite{Pandofi-Priola-Zabczyk}, \cite{Priola-Zabczyk}, and \cite{Priola-Wang}, see also our Corollary \ref{C19}. 

\begin{definition}
{\rm Let  $(P_s)$ be a transition semigroup on $B_b(E)$. A function $\psi \colon E\mapsto \mathbb{R}$ is called \emph{$(P_s)$-harmonic} if $P_s\psi =\psi$ for any $s\ge 0$.   We say that  $(P_s)$ has the \emph{Liouville  type property} if any bounded $(P_s)$-harmonic function is constant.}
\end{definition}
\subsection{$n$-fold symmetric It\^o's integral} Assume that the deterministic system \eqref{E12} is null-controllable at time $T$. Let us fix $t\ge T$ and a control operator $\mathcal{U}(t)$.   Given $y$ and a sequence $(y_l)$ of elements of $E$,  and $s\in [0,t]$, write 
\begin{equation}\label{E13}
\begin{aligned}
I_s(y) &= \int_0^s \langle u(t,y; r),   \d W(r)\rangle_H, \\
J_s(y_l,y_j) &=  \int_0^s \langle u(t,y_l;r), u(t,y_j;r)\rangle_H \d r,\\
 u(t,y;r)&= \left(\mathcal{U}(t)y\right)(r).
\end{aligned}
\end{equation}
For $n\ge 1$ we define $n$-folde symmetric integrals $I^n_s(y_1, \ldots, y_n)$ putting  $I^1_s(y)=I_s(y)$ and 
\begin{align*}
&I^{n}_s(y_{1},\ldots, y_{{n}})\\
&= \sum_{j=1}^n \int_0^s I^{n-1}_r(y_1,\ldots, y_{j-1}, y_{j+1}, \ldots y_n) \langle u(t,y_{j}; r),   \d W(r)\rangle_H. 
\end{align*}

\subsection{General result}
Given $u\in L^2(0,t;H)$ define  
$$
M_{u}(t):=\exp\left\{\int_0^t \langle u(s),\d W(s)\rangle_H -\frac 12  \int_0^t \| u(s)\|_H^2\d s\right\}.
$$

The following result is a generalization of Theorem 9.26 from \cite{DaPrato-Zabczyk2}. We denote by $D$ the operator of Fr\'echet derivative in space variable.  
\begin{theorem}\label{T17}
Assume that the deterministic system \eqref{E12} is null-controllable at time $T$.  Then for  $\phi\in B_b(E)$, $t\ge T$ and $x\in E$ we have 
\begin{equation}\label{E14}
P_t\phi(x)= \mathbb{E}\, \phi(X^0(t))M_{-u(t,x,\cdot)}(t). 
\end{equation}
Moreover, for any $n$ and $\phi\in B_b(E)$, $P_t\phi$ is n-times differentiable, and
\begin{equation}\label{E15}
\begin{aligned}
D^n  P_t\phi(x)[y_1,\ldots,y_n]&= \mathbb{E}\,\phi(X^x(t))(-1)^n I^n_t(y_1,\ldots,y_n). 
\end{aligned}
\end{equation}
Finally, we have  
$$
\left\vert  D^n P_t\phi(x)[y_1,\ldots,y_n]\right\vert \le \left\|\mathcal{U}(t)\right\|^{n}_{L(E,L^2(0,t;H))} \|\phi\|_{B_b(E)} \prod_{l=1}^n \| y_l\|_E. 
$$
\end{theorem}
\begin{proof} Let $t\ge T$. By the Girsanov theorem  
$$
W^*(s)= W(s)-\int_0^s  u(t,x;r) \d r, \qquad s\in [0,t], 
$$
is a  Wiener process under the probability measure 
$$
\d \mathbb{P}^*= M_{u(t,x;\cdot)}(t)\d \mathbb{P} = \e^{I_t(x)-\frac{1}{2}J_t(x,x)}. 
$$
Then
$$
\mathbb{E}\,\phi(X^x(t))=  \mathbb{E}^*\phi(X^x(t)) \left(M_{u(t,x;\cdot))}(t)\right)^{-1},
$$
\begin{align*}
&\left(M_{u(t,x;\cdot)}(t)\right)^{-1}\\
&= \exp\left\{ - \int_0^t  \langle u(t,x;s), \d W(s)\rangle_H + \frac12  \int_0^t \|u(t,x;s)\|_H ^2\d s\right\}\\
&=\exp\left\{ - \int_0^t  \langle u(t,x;s), \d W^*(s)\rangle_H - \frac12 \int_0^t \| u(t,x;s)\|_H ^2\d s \right\},
\end{align*}
and from  the null-controllability assumption 
\begin{align*}
X^x(t)=  \int_0^t S(t-s)a\d s+ \int_0^ t S(t-s)B \d W^*(s).
\end{align*}
Hence we have \eqref{E14}, and   moreover, 
$$
D^nP_t\phi(x)[y_1,\ldots, y_n] = \mathbb{E}\, \phi(X^0(t))D^nM_{-u(t,x,\cdot)}(t)[y_1,\ldots, y_n]. 
$$
Under the measure $\d \mathbb{P^*}= M_{-u(t,x;\cdot)}(t)\d \mathbb{P}$ the process
$$
W^*(s)= W(s)+ \int_0^s  u(t,x;r)\d r 
$$
is  a Wiener process. Moreover, 
\begin{align*}
X^0(t)&=   \int_0^t S(t-s)a\d s+ \int_0^ t S(t-s)B\left( \d W^*(s)- u(t,x;s)\d s\right)\\
&= S(t)x + \int_0^t S(t-s)a\d s+   \int_0^ t  S(t-s)B \d W^*(s). 
\end{align*}
Therefore 
$$
D^nP_t\phi(x)[y_1,\ldots, y_n]= \mathbb{E}\, \phi(X^x(t))Y(t,y_1,\ldots,y_n),
$$
where $Y(t,y_1,\ldots,y_n)$ can be obtain from  
\begin{align*}
&\frac{D^nM_{-u(t,x,\cdot)}(t)[y_1,\ldots, y_n]}{M_{-u(t,x,\cdot)}(t)}\\
& =\e^{I_t(x)+ \frac{1}{2} J_t(x,x)} D^n \e^{-I_t(x)- \frac{1}{2} J_t(x,x)}[y_1,\ldots,y_n]
\end{align*}
by replacing above $\d W(s)$ by $\d W(s)- u(t,x;s)\d s$, that is each $I_t(y)$ by $I_t(y)- J_t(y,x)$. In particular, we have  
\begin{align*}
Y(t,y)&= -\int_0^ t \langle u(t,y;s), \d W(s)- u(t,x;s)\d s\rangle_H \\
&\quad - \int_0^t\langle u(t,x;s),u(t,y;s)\rangle _H\d s\\
&=-\int_0^t \langle u(t,y;s), \d W(s)\rangle _H= -I_t(y),  
\end{align*}
and 
\begin{align*}
Y(t,y_1,y_2)&= I_t(y_1)I_t(y_2)-J_t(y_1,y_2).   
\end{align*}

In order to see that for any $n$, $Y(t,y_1,\ldots y_n)= (-1)^nI^n_t(y_1,\ldots,y_n)$, let us fix   $t$, $x$, $n$, $y_1,\ldots,y_n$, and let us write 
\begin{align*} 
Z_n(s)&:= D^n M_{-u(t,x,\cdot)}(t)[y_1,\ldots,y_n],\\
Z_n^j(s)&:= D^{n-1} M_{-u(t,x,\cdot)}(t)[y_1,\ldots,y_{j-1},y_{j+1},\ldots, y_n].
\end{align*}
Then 
\begin{align*}
\d M_{-u(t,x,\cdot)}(s)&= -M_{-u(t,x,\cdot)}(s)\left\langle u(t,x;s),\d W(s)\right\rangle _H,\\
\d Z_n(s)&= -Z_n(s)\left\langle u(t,x;s),\d W(s)\right\rangle _H \\
&\qquad - \sum_{j=1}^n Z_n^j (s)\left\langle u(t,y_j;s),\d W(s)\right\rangle _H, \\
\d \frac{1}{M_{-u(t,x,\cdot)}(t)}&= \frac{1}{M_{-u(t,x,\cdot)}(t)}\left[ \langle u(t,x;s),\d W(s)\rangle _H+\| u(t,x;s)\|^2_H\right] \d s,
\end{align*}
and consequently, for 
$$
V_n(s):= \frac{Z_n(s)}{M_{-u(t,x,\cdot)}(t)}\qquad \text{and}\qquad V_n^j(s):= \frac{Z_n^j(s)}{M_{-u(t,x,\cdot)}(t)}
$$
we have 
\begin{align*}
\d V_n(s)&= \frac{\d Z_n(s)}{M_{-u(t,x,\cdot)}(t)} + Z_n(s)\d \frac{1}{M_{-u(t,x,\cdot)}(t)}\\
&\quad  - V_n(s)\| u(t,x;s)\|_H^2 \d s -\sum_{j=1}^n V_n^j(s)\langle u(t,x;s), u(t,y_j;s)\rangle_{H}\d s\\
&= -V_n(s)\left\langle u(t,x;s),\d W(s)\right\rangle _H - \sum_{j=1}^n V_n^j (s)\left\langle u(t,y_j;s),\d W(s)\right\rangle _H\\
&\quad  + V_n(s)\left\langle u(t,x;s),\d W(s)\right\rangle _H+V_n(s)\| u(t,x;s)\|^2_H \d s\\
&\quad - V_n(s)\| u(t,x;s)\|_H^2 \d s -\sum_{j=1}^n V_n^j(s)\langle u(t,x;s), u(t,y_j;s)\rangle_{H}\d s\\
&= -\sum_{j=1}^n V_n^j (s)\left[ \left\langle u(t,y_j;s),\d W(s)\right\rangle _H+\langle u(t,x;s), u(t,y_j;s)\rangle_{H}\d s \right]. 
\end{align*}
Now, let $Y_n(s)$ and $Y_n^j(s)$ be obtained from $V_n(s)$ and $V_n^j(s)$ by replacing $\d W(s)$ by $\d W(s)- u(t,x;s)\d s$. We have  
\begin{align*}
\d Y_n(s)&=  -\sum_{j=1}^n Y_n^j (s) \left\langle u(t,y_j;s),\d W(s)\right\rangle _H, 
\end{align*}
as required. 
\end{proof}

\begin{corollary}\label{C18}
If the deterministic system is  null-controllable with vanishing energy then for any $\phi \in B_b(E)$, $n$, and $y_1,\ldots, y_n\in E$, 
$$
\sup_{x\in E}  \left \vert D^n P_t\phi(x)[ y_1,\ldots, y_n]\right\vert \to 0\quad \text{as $t\to +\infty$.}
$$
\end{corollary}

\begin{corollary}\label{C19}
If the deterministic system is null-controllable with vanishing energy then $(P_t)$ has the Liouville  property. 
\end{corollary}
\begin{proof}  Let $\psi$ be a  bounded harmonic function. Then  $P_t\psi =\psi$ for any $t\ge T$.  If the system is null-controllable with vanishing energy, then by Corollary \ref{C18},  for all $x,y\in E$ we have $DP_t\psi(x)[y]\to 0$ as $t\to +\infty$.  Since $P_t\psi=\psi$, we have $D\psi(x)[y]=0$ for $x,y\in E$, and therefore $\psi$ is constant. 
\end{proof}

\begin{corollary}\label{C110}
Under assumptions of Theorem \ref{T17}, for any $t\ge T$, the laws $\mathcal{L}(X^x(t))$, $x\in E$, are equivalent. 
\end{corollary}
\begin{proof} Taking in \eqref{E14}, as $\phi$ the  characteristic function of a Borel subset $\Gamma$ of $E$, we obtain 
$$
\mathbb{P}\left\{ X^x(t)\in \Gamma\right\}= \int_{\{X^0(t)\in \Gamma\}}M_{-u(t,x;\cdot)}(t)\d \mathbb{P}. 
$$ 
Therefore  $\mathcal{L}(X^x(t))$ is absolutely continuous with respect to $\mathcal{L}(X^0(t))$ for any $x$. We can now repeat the arguments leading to \eqref{E14}. Namely, 
\begin{align*}
X^x(t)&= \int_0^t S(t-s)a \d s+ \sum_{k} \int_0^t S(t-s)Be_k  \left[ \d W_k(s)-u_k(t, x;s)\d s\right] \\
&= \int_0^t S(t-s)a \d s + \sum_{k} \int_0^t S(t-s)Be_k   \d W^*_k(s),
\end{align*}
where 
$$
W^*_k(s)= W(s)-\int_0^s u_k(t,-x;r)\d r
$$
is a Wiener process with respect to  $\d \mathbb{P}^*=M_{u(t,x;\cdot)}(t) \d \mathbb{P}$.
Therefore 
\begin{align*}
\int_{\{ X^x(t)\in \Gamma \}}M_{u(t,x;\cdot)}(t)\d \mathbb{P}= \mathbb{P}^*\left\{ X^x(t)\in \Gamma\right\}= \mathbb{P}\left\{ X^0(t)\in \Gamma\right\}, 
\end{align*}
hence $\mathcal{L}(X^0(t))$ is absolutely continuous with respect to $\mathcal{L}(X^x(t))$. 
\end{proof}
\begin{corollary} \label{C111}
Hypothesis \ref{H11} and null-controllability at $T$ of \eqref{E12} imply that the operators $S(t)$, $t\ge T$,  on $E$ are  radonifying, and in particular compact. 
\end{corollary}
\begin{proof} By Corollary \ref{C110}, the gaussian laws $\mathcal{L}(X^x(t))$, $x\in E$, are equivalent. Let $\mathcal{H}(t)$ be the Reproducing Kernel Hilbert Space of $\mathcal{L}(X^0(t))$. By,  the Feldman--Hajek theorem, $S(t)x\in \mathcal{H}(t)$. Since the embedding $\mathcal{H}(t)\hookrightarrow E$ is radonifying, see e.g. \cite{Kuo}, the result follows. 
\end{proof}
\begin{example} In the case of the classical Liouville theorem for the Laplace operator $\Delta$ on $\mathbb{R}^d$,  $A\equiv 0$, $B\equiv \frac 12 I\in M(d\times d)$, and $W$ is a standard Wiener process on $\mathbb{R}^d$. The corresponding deterministic system is of  the form 
$$
\d Y(s)= \frac12 u(s)\d s, \qquad Y(0)=x.
$$
Clearly, it is null-controllable with vanishing energy at any time $T$ and $\mathcal{U}(t)x\equiv -\frac{2 x}{t}$. Consequently, the classical Liouville theorem saying that any bounded harmonic function on $\mathbb{R}^d$ is constant follows. 
\end{example}
\begin{remark}\label{R111}
{\rm  A natural question is if in Corollary \ref{C19} one can replace the assumption of the boundedness of a harmonic function by its appropriate integrability. Clearly, under hypothesis of Theorem \ref{T17}, for any $p,q>1$ such that $1/p+1/q=1$, there exists a constant $c(q)$ depending only on $q$ such that  
\begin{align*}
\left\vert \nabla P_t \psi(x)[y]\right\vert &=  \left\vert \mathbb{E}\, \psi(X^x(t))I_t(y)\right\vert \le \left(\mathbb{E}\left\vert \psi(X^x(t))\right\vert ^{p}\right)^{1/p}\left( \mathbb{E} \left\vert I_t(y)\right\vert ^q\right)^{1/q}\\
&\le c(q) \left(\mathbb{E}\left\vert \psi(X^x(t))\right\vert ^{p}\right)^{1/p} \left\| \mathcal{U}(t)y\right\| _{L^2(0,t;E)}. 
\end{align*}
If $\psi$ is harmonic, than in order to show that $\psi$ is constant, or equivalently, that $\nabla \psi(x)\equiv 0$, we need to show  that 
\begin{equation}\label{E19}
\lim_{t\to +\infty} \mathbb{E}\left\vert \psi(X^x(t))\right\vert ^{p} \left\| \mathcal{U}(t)y\right\| _{L^2(0,t;E)}=0.
\end{equation}
If the deterministic system is null-controllable with vanishing energy, then $\left\| \mathcal{U}(t)y\right\| _{L^2(0,t;E)}\to 0$ as $t\to +\infty$. Therefore \eqref{E19} holds if 
\begin{equation}\label{E110}
\sup _{t>0} \mathbb{E}\left\vert \psi(X^x(t))\right\vert ^{p} <+\infty.
\end{equation}
In the Wiener case, Lebesgue measure is invariant, but there are no finite invariant measures. Therefore 
$$
\mathbb{E}\left\vert \psi(X^x(t))\right\vert ^{p}\to \int_{\mathbb{R}^d} \left \vert \psi(z)\right\vert ^p \d z. 
$$
Consequently, any harmonic integrable function is constant. However, since integrable continuous functions are bounded this is not an extension of the classical property.  On the other hand, if there exists a unique finite invariant measure, say $\mu$, that is if 
$$
\int_0^t S(t-s)B\d W(s)
$$ 
converges weakly to $\mu$  as $t\to +\infty$, then the Liouville property can be stated as follows: if the deterministic system is null-controllable with vanishing energy than any harmonic function $\psi$ such that $\psi\in L^p(E, \mu)$ for a certain $p>1$ is constant. 
}
\end{remark}
\begin{remark}
{\rm If $E$ is a Hilbert space, then, see e.g. \cite{DaPrato-Zabczyk2, DaPrato-Zabczyk3},  there exists a finite invariant measure $\mu$ (say) for \eqref{E11} if and only if 
$$
\int_0^{+\infty} \|S(t)B\|^2_{L_{(HS)}(H,E)}\d t <+\infty.
$$  
If there exists the spectral gap for the generator $L$ of the transition semigroup $(P_s)$ in $L^2(E,\mu)$ then the Liouville property holds for functions in $L^p(E,\mu)$ for every $p\in(1,+\infty)$, see \cite{ania,ou-symm} for more explict conditions. 
}
\end{remark}
\section{Gradient estimate implies null-controllability}

Assume that for  all $t\ge T$ and $x\in E$  there is a integrable $Y(t,x)\colon \Omega\mapsto E^*$ such that for any $\phi\in B_b(E)$ and any $z\in E$, 
$$
DP_t\phi(x)[z]= \mathbb{E}\, \phi(X^x(t))Y(t,x)[z]. 
$$
Taking non-zero but constant $\phi$ we see that  $\mathbb{E} \, Y(t,x)[z]=0$. 

\begin{proposition} 
Assume that  $Y(t,x)[z]$ are  square integrable. Then the deterministic system is null--controllable at time $T$. Moreover, for any $x\in E$,  $\mathcal{U}(t)z= -\mathbb{E}\, V(t,x;\cdot)[z]$, $z\in E$, where $V$ is such that 
$$
\mathbb{E}\left( Y(t,x)[z]\vert\sigma(W(s)\colon s\le t)\right)= \int_0^t \langle V(t,x;s)[z], \d W(s)\rangle_H.
$$
\end{proposition}
\begin{proof}
In the first step note that the gradient formula holds for $\phi(x)=e[x]$ where $e$ is an arbitrary element of $E^*$, that is 
$$
D P_t\, e[X^x(t)][z]= \mathbb{E}\, e[X^x(t)] Y(t,x)[z]. 
$$
Let $(e_k)$ be an orthonormal basis of $H$. Given $e\in E^*$ let $\phi(x)=e[x]$. We have 
\begin{align*}
D P_t\phi(x)[z]&= D e\left[ S(t)x\right][z]= e\left[ S(t)z\right]\\
&= \mathbb{E}\, e\left[ S(t)x + \int_0^t S(t-s)B\d W(s)\right] \int_0^t \langle V(t,x;s)[z], \d W(s)\rangle_H\\
&= \mathbb{E} \sum_{k} \int_0^t e\left[ S(t-s)Be_k\right]\d W_k(s)\int_0^t \langle V(t,x;s)[z],e_k\rangle_H \d W_k(s)\\
&= \sum_{k} \int_0^t e\left[ S(t-s)Be_k\right] \langle \mathbb{E}\, V(t,x;s)[z],e_k\rangle _H\d s. 
\end{align*} 
Set $v(t,x;s)[z] = \mathbb{E} \, V(t,x;s)[z]$. Then 
$$
e\left[ S(t)z\right] = \sum_{k} \int_0^t e\left[ S(t-s)Be_k\right] \langle v(t,x;s)[z],e_k\rangle _H\d s. 
$$
Since the above identity holds for any $e\in E^*$, and 
$$
v(t,x;s)[z]= \sum_{k} \langle v(t,x;s)[z],e_k\rangle_H e_k, 
$$
we have 
\begin{align*}
S(t)z&= \sum_{k} \int_0^t  S(t-s)Be_k \langle v(t,x;s)[z],e_k\rangle _H\d s\\
&= \int_0^t S(t-s)B v(t,x;s)[z]\d s. 
\end{align*}
\end{proof}

\section{Conditions for null-controllability}
Assume that $E$ is reflexive, which is true in the case of $L^p$ spaces. Let 
$$
\mathcal{Q}(t)\colon L^2(0,t;H)\ni u \mapsto \sum_{k} \int_0^t S(t-s)Be_k   \langle u(s),e_k\rangle_H  \d s\in E. 
$$
Then the deterministic system \eqref{E12} is null-controllable at time $T$ if and only if  $\text{range}\, S(t)\subset \text{range}\, \mathcal{Q}(t)$ for $t\ge T$.

A simple proof of the following useful observation is left to the reader.   
\begin{lemma}
For any $0\le t\le t'$, $\text{range}\, \mathcal{Q}(t)\subset \text{range}\, \mathcal{Q}(t')$. Consequently, if for given $T_0\ge 0$, 
\begin{equation}\label{E31}
\textrm{range}\, S(t+T_0)\subset \textrm{range}\, \mathcal{Q}(t)\quad \textrm{for $t>0$,}
\end{equation}
then for any $T>T_0$, the deterministic system \eqref{E12} is null-controllable at time $T$. 
\end{lemma}

The following lemma  follows directly from more general result from \cite{DaPrato-Zabczyk2,DaPrato-Zabczyk3}, 
\begin{lemma}
The following conditions are equivalent:
$$
\exists \, C\colon\quad  \{ S(t)x\colon \|x\|_E\le 1\} \subset \{ \mathcal{Q}(t)u\colon \|u\|_{L^2(0,t;H)}\le C\},
$$
and 
\begin{equation}\label{E32}
\exists \, C\colon \quad  \|S^*(t)z\|_{E^*}\le C\| \mathcal{Q}(t)^*z\| _{L^2(0,t;H)}\ \ \text{for $z\in E^*$.}
\end{equation}
\end{lemma} 

\begin{corollary}System \eqref{E12} is null-controllable  at time $T$  if and only if \eqref{E32} holds true for $t\ge T$. Moreover,  $\|\mathcal{U}(t)\|_{L(E,L^2(0,t;H))} \le C$, where $C$ is a constant in \eqref{E32}.  If for $t>0$,
\begin{equation}\label{E33}
\exists \, C\colon \quad  \|S^*(t+T_0)z\|_{E^*}\le C\| \mathcal{Q}(t)^*z\| _{L^2(0,t;H)}\ \ \text{for $z\in E^*$,}
\end{equation}
then \eqref{E12} is null-controllable at any time $T>T_0$. 
\end{corollary}

Note that 
$$
\left[\mathcal{Q}(t)^*z\right](s)= \sum_k \left\langle \left[ S(t-s)B\right]^* z,e_k\right\rangle _He_k. 
$$
Therefore \eqref{E32} and \eqref{E33} have  the forms 
\begin{equation}\label{E34}
\|S^*(t)z\|_{E^*} ^2\le C^2 \sum_{k} \int_0^t  \left\langle \left[ S(s)B\right]^* z,e_k\right\rangle _H^2 \d s, \quad z\in E^*,
\end{equation}
and
\begin{equation}\label{E35}
\|S^*(t+T_0)z\|_{E^*} ^2\le C^2 \sum_{k} \int_0^t  \left\langle \left[ S(s)B\right]^* z,e_k\right\rangle _H^2 \d s, \quad z\in E^*,
\end{equation} 

\section{The case of invertible $B$}
In this section we apply our results to a well understood case, when when $E=V=H$ is a Hilbert space and $B$ is a bounded invertible linear operator on $E$, see \cite{dzsf}, \cite{DaPrato-Zabczyk2} and references therein. In this case Hypothesis  \ref{H11} and the null-controllability of \eqref{E12} can be easily verified.  Then, see Remark \ref{R12}, Hypothesis \ref{H11} holds if and only if 
\begin{equation}\label{E41}
\int_0^t \|S(s)\|^2_{(HS)}\d s <+\infty, \qquad \forall \, t>0. 
\end{equation}
Clearly the system \eqref{E12} is null-controllable at any time $T>0$, as one can takie 
$$
\left(\mathcal{U}(t)x\right)(s)= u(t,x;s)=  -\frac1 t  B^{-1}S(s)x. 
$$
For, we have 
\begin{align*}
\int_0^t S(t-s)B \frac{1}{t} B^{-1}S(s)x\d s &= \frac{1}{t} \int_0^t S(t)x\d s = S(t)x. 
\end{align*}
Note that $\eqref{E12}$ is null-controllable with uniform vanishing energy if
\begin{equation}\label{E42}
\sup_{t\ge 0}\|S(t)\|_{L(E,E)}<+\infty. 
\end{equation}
Finally, note that there is a finite invariant measure for \eqref{E11} if $S$  is exponentially stable (and in order to have a solution in $E$,  \eqref{E41} holds true). 

Consequently, by Theorem \ref{T17}  we have the following result. The theorem provides a formula for the gradient of the transition semigroup corresponding to linear equation. It is a special case of the Bismut--Elworthy--Li formula valid for non-linear diffusions, see e.g. \cite{Peszat-Zabczyk}.  
\begin{theorem}\label{T41}
Let $\phi\in B_b(E)$, $t>0$ and $x\in E$. Then  $P_t\phi$ has derivatives of all orders  and for any $y, z\in E$, 
$$
\begin{aligned}
D  P_t\phi(x)[y] &= \mathbb{E}\left[ \phi(X^x(t))\left(I_t(u)\right)\right],\\
D^2  P_t\phi(x)[y_1,y_2] &= \mathbb{E}\left[  \phi(X^x(t))\left( I_t(y_1)I_t(y_2)- J_t(y_1,y_2)\right)\right],
\end{aligned}
$$
where 
\begin{align*}
I_t(y)&= -\frac{1}{t} \sum_{k} \int_0^t \langle B^{-1} S(s)y,e_k\rangle_E \d W_k(s),\\
J_t(y_1,y_2)&= \frac{1}{t^2} \sum_k  \int_0^t \langle B^{-1} S(s)y_1, B^{-1} S(s)y_2\rangle_E \d s.
\end{align*}
Finally we have 
$$
\sup_{\phi\in B_b(E)\colon \|\phi\|_{B_b(E)}\le 1} \sup_{x\in E} \|D^n  P_t\phi(x)\|_{L(E\times E\ldots \times E)} =\textrm{O}\left(t^{-\frac{n}{2}}\right) \qquad \text{as $t\downarrow 0$.} 
$$
\end{theorem}

\begin{example}
Consider the one-dimensional Ornstein--Uhlenbeck equation
$$
\d X= \left(-\gamma  X +a\right)\d t + b\d W,
$$
where $\gamma >0$, $a\in \mathbb{R}$, and $b\not =0$. Then there is a   unique invariant measure  $\mu= \mathcal{N}(a/\gamma, b^2/(2\gamma))$. Moreover, the corresponding deterministic system is null-controllable with vanishing energy. Therefore, see Remark \ref{R111},  the following Liouville property holds true: if $\psi\in C^2$ is such that 
$$
\frac{b^2}{2}\psi''(x)+\left(-\gamma x + a\right)\psi'(x)=0, \qquad x\in \mathbb{R}, 
$$
and for a certain $p>1$, 
$$
\int_{\mathbb{R}} \vert \psi(x)\vert ^p \e^{- \gamma (x-\frac{a}{\gamma})^2/b^2}\d x<+\infty,
$$
then $\psi $ is constant. 
\end{example}
\subsection{Analytic semigroup on a Hilbert space and degenerate noise}\label{S31}
Assume that $E$ is a Hilbert space, the semigroup $S$ is analytic, and  $\lambda $ is  in the resolvent set of its generator $(A,\textrm{Dom}\,A)$. Then, the fractional powers $(\lambda-A)^{r}$, $r\in \mathbb{R}$, are well defined. Given $r>0$, let us equip 
$$
H^r:= (\lambda-A)^{-\frac{r}{2}}E = \text{Dom}\left((\lambda-A)^{\frac{r}{2}}\right), 
$$
with the scalar product and norm inherited form $E$ by the mapping $(\lambda-A)^{-\frac{r}{2}}$. 

Let $r>0$ and let $W$ be a cylindrical Wiener process in $H= H^r$, and let $B\colon H^r\to H^r$ be invertible. Consider the following linear equation 
\begin{equation}\label{E43}
\d X= \left( AX+ a\right)\d t + \d W= \left(AX+a\right)\d t + \sum_{k} Be_k \d W_k
\end{equation}
Given $t>0$ and $x\in E$ write 
\begin{equation}\label{E44}
u(t,x;s)= \begin{cases}
0&\text{for $s<\frac{t}{2}$,}\\
-\frac{2}{t} B^{-1}S(s)x&\text{for $s\ge \frac{t}{2}$. }
\end{cases} 
\end{equation}
\begin{theorem}\label{T43}
Assume that $a\in E$, and that for some or equivalently for any $t>0$, 
\begin{equation}\label{E45}
\int_0^t \left \|S(s)(\lambda -A)^{-\frac{r}{2}}\right \|^2_{L_{(\textrm{HS})}(E,E)}\d s <+\infty. 
\end{equation}
Then \eqref{E43} defines Markov family $(X^x\colon x\in E)$ on $E$. Moreover, for the corresponding  semigroup $(P_t)$; 
$$
DP_t\phi(x)[y] = \frac{2}{t}\mathbb{E}\left\{  \phi(X^x(t)) \int_{\frac{t}{2}}^t \langle B^{-1}S(s)y,  \d W(s)\rangle _{H^r}\right\}
$$
and there is a constant $C$ such that for all $y\in E$ and $t\le 1$,
\begin{equation}\label{E46}
\begin{aligned}
\mathbb{E} \left\vert \int_{\frac{t}{2}}^t \langle S(s)y,  \d W(s)\rangle _{H^r}\right\vert &= \mathbb{E}\left\vert  \sum_{k} \int_{\frac{t}{2}}^t \langle B^{-1}S(s)y, e_k\rangle_{H^r}  \d W_k(s)\right\vert \\
&\le C\|y\|_{E} t^{\frac{1-r}{2}}. 
\end{aligned}
\end{equation}
Eventually, 
$$
\left\vert D P_t\phi(x)[y]\right\vert \le C' t^{-\frac{1+r}{2}} \, \|y\|_{E}\|\phi\|_{B_b(E)}.  
$$
Finally, \eqref{E42} guarantees that the corresponding deterministic problem is null-controllable with uniform vanishing energy. 
\end{theorem}
\begin{proof} \eqref{E45} is an if and only if condition for the existence of a solution $X^x$ to \eqref{E43} for some or equivalently any initial value $x\in E$. Therefore Hypothesis \ref{H11} is fulfilled. 

We will show that the deterministic system \eqref{E12} is null-controllable at any time.  To do this, $t>0$,  and $x\in  E$. Then 
$$
S(t)x = -\int_0^t S(t-s)Bu(t,x;s) \d s,
$$
where $u(t,x;\cdot)$ is given by\eqref{E44}.  We need to show that $u(t,x;\cdot)\in L^2(0,t;H^r)$. Since $S$ is analytic $S(s)x= S(s-t/2)S(t/2)x\in H^r$ for $s\ge \frac{t}{2}$, and there is a constant $C$ such that 
$$
\|S(s)x\|_{H^r}\le C t^{-\frac{r}{2}}\|x\|_E, \qquad \frac{t}{2}\le s \le t\le 1.  
$$
In the same way we show \eqref{E46}. 
\end{proof} 
\section{Examples}
\begin{example}
Consider the following stochastic heat equation on $\mathbb{R}$,
$$
\d X= \Delta X\d t + \d W,
$$
driven by space-time white noise $W$. Thus $H=L^2(\mathbb{R})$. Since 
$$
\int_0^t S(t-s)\d W(s)
$$
is stationary in space, it cannot live in $L^2(\mathbb{R})$. It is well-defined, however in weighted space $L^2(\mathbb{R}, \rho(x)\d x)$, where 
$\rho(x)=\left(1+ x^2\right)^{-1}$. Therefore in our case $V=E= L^2(\mathbb{R}, \rho(x)\d x)$, $H=L^2(\mathbb{R})$, and $B$ is the embedding operator. It is known, see eg. \cite{Peszat-Zabczyk1}, that $\Delta$ generates $C_0$-semigroup on $L^2(\mathbb{R}, \rho(x)\d x)$, and Hypothesis \ref{H11} holds true. Unfortunately, the null-controllability hypothesis is not satisfied. For, assume that the initial condition $x$ is the constant $1$ function. Then $S(t)x=x\not \in L^2(\mathbb{R})$. On the other hand, for any $u\in L^2(0,t;L^2(\mathbb{R}))$, 
$$
\int_0^t S(t-s) u(s)\d s\in L^2(\mathbb{R}).
$$   
\end{example}
\begin{example}
Consider the infinite system of independent scalar linear equations
$$
\d X_j= -\alpha _j X_j\d t + \sigma_j \d W_j, 
$$
where $(\alpha _j)$ and $(\sigma_j)$ are sequences of strictly positive numbers. Assume that $\sup_j \sigma _j <+\infty$. Then $B\colon l^2\to l^2$, $B(x_j)= (\sigma_jx_j)$
 is a bounded linear operator. Moreover, the diagonal operator $A(x_j)= (-\alpha _jx_j)$ generates $C_0$-semigroup $S$ on $l^2$, and $S(s)(x_j)= (\e^{-\alpha _j s}x_j)$. Note that 
$$
\sum_j \mathbb{E} \left[ \int_0^t \e^{-\alpha _j (t-s)} \sigma_j \d W_j(s)\right] ^2 = \sum_j \int_0^t \e^{-2\alpha_j s} \sigma_j^2 \d s. 
$$
Therefore Hypothesis \ref{H11} holds true in the case of  $V=E=H=l^2$, provided that 
$$
\sum_{j} \frac{\sigma_j^2}{\alpha _j} <+\infty. 
$$
Let us check the null-controllability  hypothesis. Given $x=(x_j)\in l^2$ and $t>0$ we are looking for $u=(u_j)\in L^2(0,t;l^2)$ such that 
$$
\e^{-\alpha _jt}x_j= -\sigma_j \int_0^t \e^{-\alpha _j(t-s)}u_j(s)\d s. 
$$   
Clearly, 
$$
u_j (s)= \e^{-\alpha_j s}\frac{-x_j}{t \sigma_j}
$$
solves the equation. We need to verify whether  $u=(u_j)\in L^2(0,t; l^2)$. We have 
$$
\sum_j \int_0^t u_j^2(s)\d s = \sum_j \frac{x_j^2}{2 t^2\sigma_j^2\alpha _j } \left[ 1- \e^{-2\alpha _j t}\right]. 
$$
Therefore the null--controllability hypothesis holds true if  $\inf_{j} \sigma_j^2 \alpha_j >0$.  Moreover, in this case the system is null-controllable with vanishing energy.  For more information on such systems we refer the reader to \cite{DaPrato-Zabczyk2,DaPrato-Zabczyk3}. 
\end{example}
\section{Finite dimensional degenerate linear system}
Consider now a finite dimensional case $E= V= \mathbb{R}^d$ and $H= \mathbb{R}^m$.  Then \eqref{E11} has the form 
$$
\d X= \left(AX +a\right)\d t + B\d W, \qquad X(0)=x\in \mathbb{R}^d,
$$ 
where $A$ and $B$ are   $d\times d$ and $d\times m$ matrices, $W$ is a standard Wiener process in $\mathbb{R}^m$, and $a\in \mathbb{R}^d$.  Obviously, Hypothesis \ref{H11} is satisfied. Moreover, the corresponding deterministic system  is null-controllable if and only if 
\begin{equation}\label{E51}
\text{\rm rank}\left[ B, AB, \ldots, A^{d-1}B\right]=d. 
\end{equation}
For the proof of the following result we refer the reader to the work of Seidman \cite{Seidman-3}. 
\begin{theorem}\label{T51}
Assume \eqref{E51}. Let $K$ be the minimal exponent such that 
$$
\text{\rm rank}\left[ B, AB, \ldots, A^{K}B\right]=d. 
$$
Then the deterministic system is null-controllable  at any time $T>0$, and 
$$
\|\mathcal{U}(t)\|_{L(\mathbb{R}^d, L^2(0,t;\mathbb{R}^m))}=  \text{O}\left( t^{-(K+1/2)}\right) \qquad \text{as $t\to 0$}. 
$$
\end{theorem}
As a direct consequence of Theorems \ref{T17} and \ref{T51} we have the following result. 
\begin{theorem}
Assume \eqref{E41} and let $K$ be as in Theorem \ref{T41}. Then for any $t>0$ and $\phi\in B_b(\mathbb{R}^n)$, $P_t\phi$  has derivatives of all orders and 
\begin{align*}
\sup_{x\in \mathbb{R}^d} \|D^nP_t\phi (x) \|_{\mathbb{R}^{nd}}&\le \|\phi\|_{B_b(\mathbb{R}^d)} \text{O}\left( t^{-n(K+\frac{1}{2})}\right)\quad \text{as $t\to 0$.}  
\end{align*}
\end{theorem}
\subsection{Kolmogorov's diffusion}
Consider the so-called Kolmogorov diffusion introduced in \cite{kolm1}
$$
\d X_1= X_2\d t,\qquad \d X_2= \d W.
$$
We have 
$$
\d X= AX\d t + B\d W, 
$$
where 
$$
A= \left[ \begin{array}{cc} 0&1\\0&0\end{array}\right], \qquad B= \left[ \begin{array}{c} 0\\ 1\end{array}\right].
$$
Thus 
$$
\text{\rm rank}\left[ B,AB\right] = \text{\rm rank}\left[  \begin{array}{cc} 0&1\\1&0\end{array}\right] = 2. 
$$
Therefore the rank condition holds with $K=1$, and by Theorem \ref{T51}, the deterministic system is null-controllable    and the  linear mapping $\mathcal{U}(t)\colon \mathbb{R}^2\mapsto L^2(0,t;\mathbb{R})$ is such that 
$$
\|\mathcal{U}(t)\|_{L(\mathbb{R}^2; L^2(0,t;\mathbb{R}))}= \textrm{O}(t^{-3/2})\qquad \text{as $t\downarrow 0$.} 
$$
By direct calculation on can show that 
$$
\begin{aligned}
u(t,x;s)&= a(t,x)s +b(t,x),\quad s\in [0,t],\\
a(t,x)&=  \frac{x_1}{12 t^3}-\frac{35x_2}{18t^2},\\
b(t,x)&= -\frac{1}{12} \left( \frac{x_1}{2t^2}+ \frac{x_2}{3t}\right).
\end{aligned}
$$
has the desired properties. Moreover,  the system is null-controllable with vanishing energy as for the above controls 
$$
\|\mathcal{U}(t)\|_{L(\mathbb{R}^2, L^2(0,t;\mathbb{R}^2))}= \text{O}(t^{-3/2})\qquad \text{as $t\to +\infty$.}
$$
\begin{remark}
{\rm It is easy to check that in  the case of Kolmogorov diffusion  there are no finite invariant measures. Note that in the case of Kolmogorov diffusion, $\psi\in C^2(\mathbb{R}^2)$ is harmonic if 
$$
\frac{1}{2}\frac{\partial^2 \psi}{\partial x_2^2}(x_1,x_2)+ x_2\frac{\partial \psi}{\partial x_1}(x_1,x_2)=0, \qquad (x_1,x_2)\in \mathbb{R}^2. 
$$
}
\end{remark}
\section{The case of boundary noise}
\subsection{Well posedness}
Let $\mathcal{O}$ be a bounded region in $\mathbb{R}^d$.  We assume that $\mathcal{O}$ satisfies the conditions from \cite{Goldys-Peszat}. In particular it is enough to assume that $d=1$ and $\mathcal{O}$ is an interval, or $d>1$ and $\mathcal{O}$ is a bounded region with boundary of class $C^{1+\alpha}$. Given $p\in (1,+\infty)$ and $\theta \in [0,2p-1)$ set $L^p_{\theta}:= L^p(\mathcal{O},\textrm{dist}\left( \xi, \partial \mathcal{O}\right)^\theta \d \xi)$. 

Let us denote by $S$ the semigroup on $E= L^p_{\theta}$ generated by the Laplace operator $\Delta$ on $\mathcal{O}$ with homogeneous Dirichlet boundary conditions, see \cite{Goldys-Peszat}. The semigroup is given by the Geeen kernel $G$; 
$$
S(t)x(\xi)= \int_{\mathcal{O}} G(t,\xi,\eta)x(\eta )\d \eta. 
$$

Given $\lambda >0$ let $D_\lambda$ be the Dirichlet map. Let us recall that given $\lambda\ge 0$ and a function $\gamma$ on $\partial \mathcal{O}$, $u=D_\lambda \gamma$ is,  the possibly weak, unique solution to the Poisson equation 
$$
\Delta u(\xi)=\lambda u(\xi),\quad \xi\in \mathcal{O}, \qquad u(\xi)=\gamma(\xi), \quad \xi\in \partial \mathcal{O}. 
$$

Let $H\hookrightarrow L^2(\partial \mathcal{O})$ be a Hilbert space, and let $(e_k)$ be an orthonormal basis of $H$. Finally let  $(W_k)$ be a sequence of independent real valued standard Wiener processes defined on a probability space $(\Omega,\mathfrak{F},\mathbb{P})$.  

Let us recall, see \cite{Goldys-Peszat},  that the boundary problem 
\begin{equation}\label{E71}
\begin{aligned}
\frac{\partial X}{\partial s}(s,\xi)&= \Delta X(s,\xi), \qquad s>0, \ \xi\in \mathcal{O},\\
X(0,\xi)&= x(\xi),\qquad \xi\in \mathcal{O},\\
X(s,\xi)& = \sum_{k} e_k(\xi)\frac{\d W_k}{\d s}(s), \qquad s>0, \ \xi\in \partial \mathcal{O},
\end{aligned}
\end{equation}
can be written in form \eqref{E11}. Moreover, it is well posed on the state space $L^p_{\theta}$ if and only if, for a certain  or equivalently for every $T>0$, 
\begin{equation}\label{E72}
\begin{aligned}
&\mathcal{J}_T(\{e_k\},p,\theta)\\
&= \int_{\mathcal{O}}\left[  \sum_{k} \int_0^T \left( (\lambda -A)S(s)D_\lambda e_k\right)^2(\xi) \d s \right] ^{p/2} \textrm{dist}\left(\xi,\partial \mathcal{O}\right)^\theta \d \xi
\end{aligned}
\end{equation}
is finite.  Here $A$ is the Laplace operator with homogeneous Dirichlet boundary conditions.  Moreover,  it solution $X^x$ is given by the formula 
\begin{align*}
X^x(s,\cdot)&= S(s)x(\cdot)+ \sum_k \int_0^s (\lambda -A)S(s-r)D_\lambda e_k (\cdot)\d W_k(r)\\
&=\int_{\mathcal{O}}G(s,\cdot, \eta )x(\eta )\d \eta \\
&\quad  - \sum_k \int_0^ s \int_{\partial \mathcal{O}}\frac{\partial }{\partial \mathbf{n}^a(\eta )} G(s-r, \cdot, \eta)e_k(\eta )\d \mathbf{s}(\eta)\d W_k(r),
\end{align*}
where $\mathbf{n}=\left(\mathbf{n}^1,\ldots,\mathbf{n}^d\right)$ is the outward pointing unit normal vector to the boundary $\partial \mathcal{O}$,  and $\mathrm{s}$ is the surface measure.   

For more information on the deterministic and stochastic boundary problems see e.g. \cite{bgpr,  DaPrato-Zabczyk,  Fattorini-Russell, Goldys-Peszat}.

Hypothesis \ref{H11} can be written now in the following form: 
\begin{hypothesis}\label{H71} 
For any $T>0$,  $\mathcal{J}_T(\{e_k\},p,\theta)$ given by \eqref{E72} is finite.  
\end{hypothesis}

We define the transition semigroup on $B_b(E)=B_b(L^p_{\theta})$;  $P_s\phi(x)= \mathbb{E}\, \phi \left(X^x(s)\right)$.  In the case of the boundary noise  problem the null-controllability criterion of the corresponding deterministic problem  can be formulated as follows. 
\begin{hypothesis}\label{H72}
Assume that for each $t>0$ there is a bounded linear operator $\mathcal{U}(t)\colon L^p_{\theta}\mapsto L^2(0,t;H)$, $(\mathcal{U}(t)x)(s)= u(t,x;s)$ such that for any $x\in L^p_{\theta}$, 
\begin{align*}
& \int_{\mathcal{O}}G(t,\cdot, \eta )x(\eta )\d \eta \\
&\quad  =\sum_k \int_0^ t \int_{\partial \mathcal{O}}\frac{\partial }{\partial \mathbf{n}^a(\eta)} G(t-s, \cdot, \eta)e_k(\eta )\d \mathbf{s}(\eta ) \langle u(t,x;s),e_k\rangle _H\d s.  
\end{align*}
\end{hypothesis}
\begin{remark} 
Let us note that if $\mathcal{O}$ is unbounded then it is not possible that Hypotheses \ref{H71} and \ref{H72} hold  simultaneously.  For, by Corollary \ref{C111} the semigroup $S$ need to be compact. 
\end{remark}

\begin{remark}
Note that 
$$
G(t,\xi,\eta)= \sum_{l}\e^{-\lambda_{l}t}f_l(\xi)f_l(\eta),
$$
where $(f_l)$ is the orthonormal basis of $L^2(\mathcal{O})$ of eigenvectors of the Laplace operator and $(-\lambda_l)$ are the corresponding eigenvalues. Therefore the null-controllability condition has the form 
\begin{align*}
 \langle x, f_l\rangle_{L^2(\mathcal{O})}&= \sum_{k} \int_0^t \e^{\lambda_l s} \langle u(t,x;s),e_k\rangle _H\d s \int_{\partial \mathcal{O}} \frac{\partial f_l }{\partial \mathbf{n}^a(\eta)} (\eta)e_k(\eta )\d \mathbf{s}(\eta ). 
\end{align*}
and it is enough to verify this identity for any $x\in L^2(\mathcal{O})$ and of course any $l$, see Lemmas \ref{L76} and \ref{L77} below. 
\end{remark}

\subsection{The case of white noise}
The following result is proved in \cite{Seidman-2}. Its formulation is  adapted to notations of this paper. 
\begin{theorem}\label{T75}
Then there is a bounded linear control operator $\mathcal{U}(t)$ acting from $L^2(\mathcal{O})$ into $L^2(0,t;L^2(\partial \mathcal{O}))$, and 
\begin{equation}\label{E73}
\log\, \|\mathcal{U}(t)\|_{L(L^2(\mathcal{O}),L^2(0,t;L^2(\partial \mathcal{O})))}=  \mathcal{O}( 1/t)\quad \text{as $t\downarrow 0$.} 
\end{equation}
\end{theorem} 
The following result was established in \cite{Goldys-Peszat}[Lemma 5.1].  
\begin{lemma}\label{L76}
Let $p>1$, $\theta \in [0,2p-1)$. For any $t>0$, $S(t)$ is a bounded linear operator from  $L^p_{\theta}$ into $L^p_{0}=L^p(\mathcal{O})$. Moreover, the operator norm is of order $t^{-\frac{\theta}{2p}}$, that is there exists a $C>0$ such that $\|S(t)\|_{L(L^p_{\theta},L^p_{0})}\le Ct^{-\frac{\theta}{2p}}$ for $t\le 1$.

\end{lemma}
Since $\mathcal{O}$ is bounded we have also the following result. 
\begin{lemma}\label{L77}
Let  $p>1$. Then for any $t>0$, $S(t)$ is a bounded linear operator from $L^p_0= L^p(\mathcal{O})$ into $L^2_0= L^2(\mathcal{O})$. Moreover, there exists a $C>0$ and $\alpha (p)\in [0,+\infty)$ such that $\|S(t)\|_{L(L^p_{0},L^2_{0})}\le Ct^{-\alpha(p)}$ for  $t\le 1$.
\end{lemma}

The main result of this section is the following theorem 
\begin{theorem}
(i)  Assume that $\mathcal{O}=(a,b)$ is a bounded open interval in $\mathbb{R}$, and that $W_1,W_2$ are independent Wiener processes. Let $p>1$ and $\theta \in (p-1,2p-1)$. Then the problem \eqref{E71} with boundary conditions $X(t,a)=\dot W_1(t)$ and $X(t,b)=\dot W _2(t)$ defines a Markov family  on $E= L^p_\theta$. 

\smallskip \noindent (ii) Assume that $\mathcal{O}$ is a bounded region in  $\mathbb{R}^2$ with boundary $C^{1,\alpha}$ for some $\alpha >0$. Let $(e_k)$ be an orhonormal basis of  $L^2(\partial \mathcal{O})$ and let $(W_k)$ be a sequence of independent Wiener processes. Let $p>1$ and $\theta \in \left(\frac{3p}{2}-1,2p-1\right)$. Then the problem \eqref{E71}  defines a Markov family  on $E= L^p_\theta$. 

Moreover, in both cases for any $\phi \in B_b(E)= B_b(L^p_\theta)$ and $t>0$,  $P_t\phi$ has derivatives of all orders  and 
$$
\sup_{x\in E}\|D^n P_t\phi(x)\|_{L(E\times E\times \ldots \times E)}\le C\e^{\frac{Cn}{t}}\sup_{x\in E}\|\phi\|_{B_b(E)}, \qquad \phi \in B_b(E),
$$
with $C>0$ independent of $\phi$, $t\le 1$, and $n$. 
\end{theorem}
\begin{proof} Parts $(i)$ and $(ii)$ are related to Hypothesis \ref{H11}, or Hypothesis \ref{H71},  and are know, see \cite{Goldys-Peszat}[Propositions 8.1, 8.7, 8.10],  and they are recall here only for the sake of completeness.  To show the last claim, we need to verify the  null-controllability Hypothesis \ref{H72}.  Let $x\in L^p_\theta$ and $t>0$, By Lemmas \ref{L76} and \ref{L77},   $S(t/2)x\in L^2(\mathcal{O})$ and 
$$
\|S(t/2)x\|_{L^2} \le C_1  t^{-\frac{d+ \theta}{2p}}\|x\|_{L^p_{\theta}}. 
$$
By Theorem \ref{T75} there is a 
$$
\mathcal{U}(t/2)\colon L^2(\mathcal{O}) \mapsto L^2(0,t/2, L^2(\partial \mathcal{O})
$$
such that for the control $u(s)= \mathcal{U}(t)[S(t/2)x](s)$, $s\in [0,t/2]$,   $u\in L^2(0,t/2;L^2(\partial \mathcal{O}))$, we have 
\begin{align*}
0&= S(t)x+ \int_0^{\frac{t}{2}} (\lambda -A)S(t/2-s)D_\lambda u(s)\d s\\
&= S(t)x+ \int_0^{t} (\lambda -A)S(t-s)D_\lambda T_tu(s)\d s,
\end{align*}
where $T_{t}\colon L^2(0,t/2;L^2(\partial \mathcal{O}))\mapsto L^2(0,t;L^2(\partial \mathcal{O}))$ is bounded linear operator defined as follows   
$$
T_tu(s)= \begin{cases} 0&\text{for $s\le t/2$,}\\
u(s-t/2)&\text{for $s\in [t/2,t]$.}
\end{cases}
$$
\end{proof} 
\subsection{The case of coloured noise} 
White noise case is restricted to subdomains of $\mathbb{R}^d$, $d=1,2$.  In this section we consider the case of the so-called coloured noise, that is the case where $(e_k)$ is not an orthonormal basis of $L^2(\partial \mathcal{O})$. Assume that $\mathcal{O}$ is a bounded domain in $\mathbb{R}^d$ with $C^{1,\alpha}$, $\alpha >0$,  boundary $\partial \mathcal{O}$. By  Proposition 8.10 from \cite{Goldys-Peszat}, we have the following result dealing with  the verification of Hypothesis \ref{H71}. 
\begin{proposition}\label{P79}
Assume that 
\begin{equation}\label{E74}
\sum_{k} \sup_{\eta \in \partial \mathcal {O}}e_k^2(\eta)<+\infty
\end{equation}
and $1<p<+\infty$ and $\theta \in (p-1,2p-1)$. Then the boundary problem \eqref{E71}  defines a Markov family with continuous trajectories in   $E=L^p_\theta:= L^p_{\theta,0}$. 
\end{proposition}

What is left is to find conditions on $(e_k)$ which guarantee for a given $t>0$, the null-controllability of the corresponding deterministic boundary problem 
\begin{equation}\label{E75}
\begin{aligned}
\frac{\partial Y}{\partial s}(s,\xi)&= \Delta Y(s,\xi), \qquad s\in (0,t), \ \xi\in \mathcal{O},\\
Y(0,\xi)&= x(\xi),\qquad \xi\in \mathcal{O},\\
Y(s,\xi)& = \sum_{k} e_k(\xi)u_k(t,x;s), \qquad s\in (0,t), \ \xi\in \partial \mathcal{O}. 
\end{aligned}
\end{equation}
Taking into account Lemmas \ref{L76} and \ref{L77}, it is enough to establish null-controllability in unweighted space $L^2(\mathcal{O})$.

\subsubsection{The case of  a unit ball in $\mathbb{R}^d$} Here we give an example of an equation with colored noise satisfying \eqref{E74} and such that the corresponding boundary problem \eqref{E75} is null-controllable in $L^2(\mathcal{O})$. Namely, let $\mathcal{O}=B_d(0,1)$ be the  unit ball in $\mathbb{R}^d$. Then $\partial \mathcal{O}=\mathbb{S}^{d-1}$.  In polar coordinates the eigenvectors of the Laplace operator have the form 
$$
H_{n, \mathbf{k}}(r, \theta)= h_n(r)f_{\mathbf{k}}(\theta), \qquad r\in [0,1], \ \theta \in \mathbb{S}^{d-1},
$$
$n\in \mathbb{N}$ and $\mathbf{k}\in \mathbb{N}_*^{d-1}\setminus\{0\} := (\mathbb{N}\cup\{0\})^{d-1}\setminus \{(0,\ldots, 0\}$. 

Above  $(h_n)$ are \emph{Bessel functions} and $(f_\mathbf{k})$ is an orthonormal basis of $\mathbb{S}^{d-1}$; the-so called \emph{harmonics}. For our purposes it is important that  
\begin{equation}\label{E76}
\sup_{\mathbf{k}} \sup_{\theta \in \mathbb{S}^{d-1}} \vert f_{\mathbf{k}}(\theta)\vert <+\infty. 
\end{equation}
In the theorem $\vert \mathbf{k}\vert ^2=\sum_{i=1}^{d-1}\mathbf{k}_i^2$.  
\begin{theorem}
Assume that $e_\mathbf{k}=a_\mathbf{k}f_\mathbf{k}$,  $\mathbf{k}\in \mathbb{N}_*^{d-1}\setminus\{0\}$, where $(a_{\mathbf{k}})$ is a sequence of real numbers such that 
\begin{equation}\label{E77}
\sum_{\mathbf{k}} a_\mathbf{k}^2<+\infty,
\end{equation}
and 
\begin{equation}\label{E78}
\forall\, t>0\ \exists \, \beta(t)>0\colon\ \forall\, \mathbf{k} \qquad  a_\mathbf{k}^{-2} \e^{-t\vert \mathbf{k}\vert ^2}  \le \beta(t). 
\end{equation}
Let $p>1$ and $\theta \in \left(p-1,2p-1\right)$. Then the problem \eqref{E71}  defines a Markov family  on $E= L^p_\theta$. The corresponding deterministic system is null-controllable at any time $T>0$, and for any $\phi \in B_b(E)= B_b(L^p_\theta)$ and $t>0$,  $P_t\phi$ has derivatives of all orders. Finally  
\begin{equation}\label{E79}
\sup_{x\in E}\|D^nP_t\phi(x)\|_{L(E\times \ldots\times E)}\le C(t)^n \sup_{x\in E}\|\phi\|_{B_b(E)}, \  \phi \in B_b(E),
\end{equation}
with a constant  $C(t)>0$ independent of $\phi$ and $n$. Moreover, there are a function $C_1$ and a constant $c>0$ such that  $C(t)= C_1(t/2)\beta(ct/2)$ and 
$\log C_1(t)=\mathrm{O}(1/t)$ as $t\downarrow 0$.  
\end{theorem}
\begin{proof} Assume that $p>1$ and $\theta\in (p-1,2p-1)$. Then  \eqref{E76}, \eqref{E77},  and Proposition \ref{P79} guarantee that \eqref{E71} defines Markov family in $L^p_\theta$.  By Lemmas \ref{L76} and \ref{L77} null-controllability in $L^2_0=L^2(B(0,1))$ ensures null-controllability in $E=L^p_\theta$. Therefore, we only need to verify null-controllability in $L^2(B(0,1))$ and the desired estimates  for the control operator. 

We have 
$$
G(t,r,\theta, r', \theta' )= \sum_{n,\mathbf{k}} \e^{-\lambda_{n,\mathbf{k}}t} h_n(r)f_\mathbf{k}(\theta) h_n(r')f_\mathbf{k}(\theta'), 
$$
where $(-\lambda_{n,\mathbf{k}})$ is the corresponding  sequence of  eigenvalues.  Then 
$$
S(t)x = \sum_{n,\mathbf{k}} \e^{-\lambda_{n,\mathbf{k}}t}\langle x,H_{n,\mathbf{k}}\rangle _{L^2(B(0,1))} H_{n,\mathbf{k}}. 
$$
Note that 
$$
u=\sum_\mathbf{k}  \langle u,e_\mathbf{k}\rangle_H e_\mathbf{k} = \sum_\mathbf{k} \langle u,e_\mathbf{k}\rangle _H a_\mathbf{k}f_\mathbf{k}.
$$
Thus $\langle u,f_\mathbf{k}\rangle_{L^2(\mathbb{S}^{d-1})}= \langle u,e_\mathbf{k}\rangle_H a_\mathbf{k}$,  and  $\langle u,e_\mathbf{k}\rangle_H= a_\mathbf{k}^{-1} \langle u,f_\mathbf{k}\rangle_{L^2(\mathbb{S}^{d-1})}$.

Our aim is to verify  condition  \eqref{E35}. We have 
\begin{align*}
\langle (S(s)B)^*x, e_\mathbf{k}\rangle _H&= \langle x, S(s)Be_\mathbf{k}\rangle _{L^2(\mathcal{O})}= a_\mathbf{k} \langle x, S(s)Bf_\mathbf{k}\rangle _{L^2(\mathcal{O})}\\
&= -a_\mathbf{k} \langle S(s)x, AD_0f_\mathbf{k}\rangle _{L^2(\mathcal{O})}\\
&= -a_\mathbf{k} \int_{\mathbb{S}^{d-1}}\frac{\partial }{ \partial \mathbf{n}^a(\eta)} S(s)x (\eta)f_\mathbf{k}(\eta)\d \mathbf{s}(\eta). 
\end{align*}
Since 
$$
\frac{\partial }{ \partial \mathbf{n}^a(\eta)} S(s)x = \sum_{n, \mathbf{j}} \e^{-\lambda_{n,\mathbf{j}}s} \langle x,H_{n,\mathbf{j}}\rangle_{L^2(\mathcal{O})} f_\mathbf{j} \frac{\d h_n}{\d x}(1),
$$
we have 
\begin{align*}
\sum_\mathbf{k} \langle (S(s)B)^*x, e_\mathbf{k}\rangle _H^2&=\sum_{\mathbf{k}} a_\mathbf{k}^2 \left(\sum_{n} \e^{-\lambda_{n,\mathbf{k}}s} \langle x,H_{n,\mathbf{k}}\rangle_{L^2(\mathcal{O})}  \frac{\d h_n}{\d x}(1)\right)^2.
\end{align*}
We need to show that for ant $T_0>0$ there is a function $\tilde C=\tilde C(t)$  such that   for $x\in L^2(B(0,1))$, 
\begin{align*}
&\sum_{n,k} \e^{-2\lambda_{n,\mathbf{k}}(t+T_0)}\langle x,H_{n,\mathbf{k}}\rangle_{L^2(\mathcal{O})}^2 \\
&\le \tilde C(t) \int_0^t \sum_{\mathbf{k}} a_\mathbf{k}^2 \left(\sum_{n} \e^{-\lambda_{n,\mathbf{k}}s} \langle x,H_{n,\mathbf{k}}\rangle_{L^2(\mathcal{O})}  \frac{\d h_n}{\d x}(1)\right)^2\d s\\
&= \tilde C(t)\sum_{n,m, \mathbf{k}} R_\mathbf{k}(a_\mathbf{k}, n,m,t) \e^{-\lambda_{n,\mathbf{k}}t}   \e^{-\lambda_{m,\mathbf{k}}t}\langle x,H_{n,\mathbf{k}}\rangle _{L^2(\mathcal{O})} \langle x,H_{m,\mathbf{k}}\rangle _{L^2(\mathcal{O})}, 
\end{align*}
where 
\begin{align*}
R_\mathbf{k}(a_\mathbf{k},n,m,t)&:= a_\mathbf{k}^2 \frac{1}{\lambda_{n,\mathbf{k}} +\lambda_{m,\mathbf{k}}} \left(\e^{(\lambda_{n,\mathbf{k}}+ \lambda_{m,\mathbf{k}})t}-1 \right)  \frac{\d h_n}{\d x}(1)\frac{\d h_m}{\d x}(1).
\end{align*}
Since $H_{n,\mathbf{k}}$ and $H_{m,\mathbf{j}}$ are orthogonal for $\mathbf{k}\not = \mathbf{j}$ we have to show that for any $\mathbf{k}$,  

\begin{equation}\label{E710}
\sum_{n} \e^{-\lambda_{n,\mathbf{k}}T_0} b_{n,\mathbf{k}}^2  \le \tilde C(t) \sum_{n,m} R_\mathbf{k}(a_\mathbf{k},n,m,t) b_{n,\mathbf{k}}b_{m,\mathbf{k}}, 
\end{equation}
where  $b_{n,\mathbf{k}} := \e^{-\lambda_{n,\mathbf{k}}t} \langle x,H_{n,\mathbf{k}}\rangle_{L^2(\mathcal{O})}$. By Theorem \ref{T75} there is a function $C_1$  of order  $\log C_1(t)=\mathrm{O}(1/t)$ such that   for any $x\in L^2(B(0,1))$, 
$$
\sum_{n}  b_{n,\mathbf{k}}^2  \le C_1(t) \sum_{n,m} R_\mathbf{k}(1,n,m,t) b_{n,\mathbf{k}}b_{m,\mathbf{k}} \qquad \text{for any $\mathbf{k}$}.  
$$
Hence \eqref{E710} holds true with 
$$
\tilde C(t)= C_1(t)\sup  _{n,\mathbf{k}} \e^{-\lambda_{n, \mathbf{k}}T_0} a_\mathbf{k}^{-2}. 
$$
Since $\lambda_{n,\mathbf{l}} \asymp c_1 ( n^2 + \vert \mathbf{k}\vert ^2)$, we have 
$$
\sup_{n,\mathbf{k}} \e^{-\lambda_{n, \mathbf{k}}T_0} a_\mathbf{k}^{-2}\le  \beta (cT_0),
$$
and consequently we have the desired estimate with $\tilde C(t)=C_1(t)\beta(cT_0)$.  What is left is to evaluate  $C(t)$ appearing in \eqref{E79}. We are looking for $C$ such that 
$$
\vert S(t)x\vert_{L^2(B(0,1))}^2 \le C(t) \vert \mathcal{Q}^*(t)x\vert _{L^2(0,t;H)}^2. 
$$ 
We have 
$$
\vert S(t+T_0)x\vert_{L^2(B(0,1))}^2 \le C_1(t) \beta(cT_0)  \vert \mathcal{Q}^*(t)x\vert _{L^2(0,t;H)}^2. 
$$
Since 
\begin{align*}
 \vert \mathcal{Q}^*(t)x\vert _{L^2(0,t;H)}^2= \int_0^t \left\vert  S(t-s)B x\right\vert _H\d s 
 & \le \int_0^{t+T_0} \left\vert  S(t+T_0-s)B x\right\vert _H\d s\\
 &= \vert \mathcal{Q}^*(t+T_0)x\vert _{L^2(0,t+T_0;H)}^2, 
\end{align*}
we have $C(t+ T_0)=C_1(t)\beta(cT_0)$. In particular $C(t)= C_1(t/2) \beta(ct/2)$. 
\end{proof}

\end{document}